\documentclass{llncs}
\usepackage{amsmath,amsfonts,amscd,amssymb,latexsym}

\newtheorem{algorithm}[theorem]{Algorithm}

\newcommand{\defterm}[1]{\emph{#1}}

\begin{document}

%%%%%%% Paper Info %%%%%%%
\title{Genetic Algorithms for Word Problems in Partially Commutative Groups}
\titlerunning{GAs for Word Problems in Partially Commutative Groups}
\author{Matthew J. Craven}\institute{Mathematical Sciences, University of Exeter,\\North Park Road,
Exeter EX4 4QF, UK.}

\maketitle\pagestyle{headings}
\begin{abstract}
We describe an implementation of a genetic algorithm on partially
commutative groups and apply it to the double coset search problem
on a subclass of groups. This transforms a combinatorial group
theory problem to a problem of combinatorial optimisation. We obtain
a method applicable to a wide range of problems and give results
which indicate good behaviour of the genetic algorithm, hinting at
the presence of a new deterministic solution and a framework for
further results.
\end{abstract}

\section{Introduction}
\subsection{History and Background}
Genetic algorithms (hereafter referred to as GAs) were introduced by
Holland \cite{Holland1975} and have enjoyed a recent renaissance in
many applications including engineering, scheduling and attacking
problems such as the travelling salesman and graph colouring
problems. However, the use of GAs in group theory
\cite{BBB03,Miasnikov99,2Miasnikov} has been in operation for a
comparatively short time.

This paper discusses an adaptation of GAs for word problems in
combinatorial group theory. We work inside the Vershik groups
\cite{Vershik}, a subclass of partially commutative groups (also
known as graph groups \cite{VanWyk} and trace groups). We omit a
survey of the theory of the groups here and focus on certain
applications.

There exists an explicit solution for many problems in this setting.
The biautomaticity of the partially commutative groups is
established in \cite{VanWyk}, so as a corollary the conjugacy
problem is solvable. Wrathall \cite{Wrathall1} gave a fast algorithm
for the word problem based upon restricting the problem to a monoid
generated by group generators and their formal inverses. In
\cite{Wrathall2}, an algorithm is given for the conjugacy problem;
it is linear time by a stack-based computation model.

Our work is an experimental investigation of GAs in this setting to
determine why they seem to work in certain areas of combinatorial
group theory and to determine bounds for what happens for given
problems. This is done by translating given word problems to ones of
combinatorial optimisation.

\subsection{Partially Commutative Groups and Vershik Groups}
Let $X=\{x_1,x_2,\ldots,x_n\}$ be a finite set and define the
operation of multiplication of $x_i,x_j\in X$ to be the
juxtaposition $x_ix_j$. As in \cite{Wrathall2}, we specify a
\defterm{partially commutative group} $G(X)$ by $X$ and the
collection of all elements from $X$ that \defterm{commute}; that is,
the set of all pairs $(x_i,x_j)$ such that $x_i,x_j\in X$ and
$x_ix_j=x_jx_i$. For example, take $X=\{x_1,x_2,x_3,x_4\}$ and
suppose that $x_1x_4=x_4x_1$ and $x_2x_3=x_3x_2$. Then we denote
this group $G(X)=\left<X:[x_1,x_4],[x_2,x_3]\right>$.

The elements of $X$ are called \defterm{generators} for $G(X)$. Note
that for general $G(X)$ some generators commute and some do not, and
there are no other non-trivial relations between the generators. We
concentrate on Vershik groups, a particular subclass of the above
groups. For a set $X$ with $n$ elements as above, the
\defterm{Vershik group of rank $n$} over $X$ is given by
\[
V_n= \left< X : [x_i,x_j] \;\hbox{ if }\; |i-j| \ge 2\right>.
\]
For example, in the group $V_4$ the pairs of elements that commute
with each other are $(x_1,x_3),(x_1,x_4)$ and $(x_2,x_4)$. We may
also write this as $V(X)$ assuming an arbitrary set $X$. The
elements of $V_n$ are represented by group \defterm{words} written
as products of generators. The \defterm{length}, $l(u)$, of a word
$u\in V_n$ is the minimal number of single generators from which $u$
can be written. For example $u=x_1x_2x_4\in V_4$ is a word of length
three. We use $x_i^{\mu}$ to denote $\mu$ successive multiplications
of the generator $x_i$; for example, $x_2^4=x_2x_2x_2x_2$. Denote
the empty word $\varepsilon\in V_n$.

For a subset, $Y$, of the set $X$ we say the Vershik group $V(Y)$ is
a \defterm{parabolic subgroup} of $V(X)$. It is easily observed that
any partially commutative group $G$ may be realised as a subgroup of
a Vershik group $V_n$ of sufficiently large rank $n$.

Vershik \cite{Vershik} solved the word problem in $V_n$ by means of
reducing words to their \defterm{normal form}. The Knuth-Bendix
normal form of a word $u\in V_n$ of length $l(u)$ may be thought of
as the ``shortest form'' of $u$ and is given by the unique
expression
\begin{eqnarray*} \overline{u} = x_{i_1}^{\mu_1}
x_{i_2}^{\mu_2} \ldots x_{i_k}^{\mu_k}
\end{eqnarray*}
such that all $\mu_i\neq 0,l(\overline{u})=\sum|\mu_i|$ and
\begin{enumerate}
\item[i)]  if $i_j=1$ then $i_{j+1}>1$;

\item[ii)]  if $i_j=m<n$ then $i_{j+1} = m-1$ or $i_{j+1}>m$;

\item[iii)]  if $i_j=n$ then $i_{j+1} = n-1$.
\end{enumerate}
The name of the above form follows from the Knuth-Bendix algorithm
with ordering $x_1<x_1^{-1}<x_2<x_2^{-1}<\ldots<x_n<x_n^{-1}$. We
omit further discussion of this here; the interested reader is
referred to \cite{KnuthBendix} for a description of the algorithm.

The algorithm to produce the above normal form is essentially a
restriction of the stack-based (or heap-based) algorithm of
\cite{Wrathall1} to the Vershik group, and we thus conjecture that
the normal form of a word $u\in V_n$ may be computed efficiently in
time $O\left(l(u)\log l(u)\right)$ for the ``average case''. From
now on we write $\overline{u}$ to mean the normal form of the word
$u$. For a word $u\in V_n$, we say that
\[
RF(u)=\{x_i^\alpha : l(\overline{ux_i^{-\alpha
}})=l(\overline{u})-1,\alpha =\pm 1\}
\]
is the \defterm{roof of $u$} and
\[
FL(u)=\{x_i^\alpha : l(\overline{x_i^{-\alpha
}u})=l(\overline{u})-1,\alpha =\pm 1\}
\]
is the \defterm{floor of $u$}. The roof (and floor) of $u$
correspond to the generators which may be cancelled after their
inverses are juxtaposed to the right (and left) end of $u$ to create
the word $u'$ and $u'$ is reduced to its normal form
$\overline{u'}$. For example, if $u=x_1^{-1}x_2x_6x_5^{-1}x_4x_1$
then $RF(u)=\{x_1,x_4\}$ and $FL(u)=\{x_1^{-1},x_6\}$.

\section{Statement of Problem}\label{prs}
Given a Vershik group $V_{n}$ and two words $a,b$ in the group, we
wish to determine whether $a$ and $b$ lie in the same double coset
with respect to given subgroups. In other words, consider the
following problem:

\subsubsection*{The Double Coset Search Problem
(DCSP)} Given two parabolic subgroups $V(Y)$ and $V(Z)$ of a Vershik
group $V_n$ and two words $a,b\in V_n$ such that $b\in
V(Y)\,a\,V(Z)$, find words $x\in V(Y)$ and $y\in V(Z)$ such that
$b=xay$.\newline

We attack this group-theoretic problem by transforming it into one
of combinatorial optimisation. In the following exposition, an
\defterm{instance} of the DCSP is specified by a pair $(a,b)$ of
given words, each in $V_n$, and the notation $\mathcal{M}((a,b))$
denotes the set of all \defterm{feasible solutions} to the given
instance. We will use a GA to iteratively produce ``approximations''
to solutions to the DCSP, and denote an ``approximation'' for a
solution $(x,y)\in\mathcal{M}((a,b))$ by $(\chi,\zeta)\in V(Y)\times
V(Z)$.

\subsubsection*{Combinatorial Optimisation DCSP}

\begin{itemize}
\item[] Input:\,\,\, Two words $a,b\in V_n$.
\item[] Constraints:\,\,\, $\mathcal{M}((a,b))=\{(\chi,\zeta)\in
V(Y)\times V(Z):\chi a\zeta\doteq b\}$.
\item[] Costs:\,\,\, The function $C((\chi,\zeta))=l(\overline{\chi
a\zeta b^{-1}})\geq 0$.
\item[] Goal:\,\,\, Minimise $C$.
\end{itemize}

The cost of the pair $(\chi,\zeta)$ is a non-negative integer
imposed by the above function $C$. The length function defined on
$V_n$ takes non-negative values; hence an \defterm{optimal solution}
for the instance is a pair $(\chi,\zeta)$ such that
$C((\chi,\zeta))=0$. Therefore our goal is to minimise the cost
function $C$.

As an application of our work, note that the Vershik groups are
inherently related to the braid groups, a rich source of primitives
for algebraic cryptography. In particular, the DCSP in the Vershik
groups is an analogue of an established braid group primitive. The
reader is invited to consult \cite{KoSECANTS02} for further details.

In the next section we expand these notions and detail the method we
use to solve this optimisation problem.

\section{Genetic Algorithms on Vershik Groups}
\subsection{An Introduction to the Approach}
For brevity we do not discuss the elementary concepts of GAs here,
but refer the reader to \cite{Holland1975,Michalewicz} for a
discussion of GAs and remark that we use standard terms such as
\defterm{cost-proportionate selection} and \defterm{reproductive method}
in a similar way.

We give a brief introduction to our approach. We begin with an
initial population of ``randomly generated'' pairs of words, each
pair of which is treated as an approximation to a solution
$(x,y)\in\mathcal{M}((a,b))$ of an instance $(a,b)$ of the DCSP. We
explicitly note that the GA does not know either of the words $x$ or
$y$. Each pair of words in the population is ranked according to
some cost function which measures how ``closely'' the given pair of
words approximates $(x,y)$. After that we systematically imitate
natural selection and breeding methods to produce a new population,
consisting of modified pairs of words from our initial population.
Each pair of words in this new population is then ranked as before.
We continue to iterate populations in this way to gather steadily
closer approximations to a solution $(x,y)$ until we arrive at a
solution (or otherwise).

\subsection{The Representation and Computation of Words}
We work in $V_n$ and two given parabolic subgroups $V(Y)$ and
$V(Z)$, and wish the GA to find an exact solution to a posed
problem. We naturally represent a group word $u=x_{i_1}^{\mu_1}
x_{i_2}^{\mu_2} \ldots x_{i_k}^{\mu_k}$ of arbitrary length by a
string of integers, where we consecutively map each generator of the
word $u$ as follows:
\begin{eqnarray*}
x_i^{\epsilon_i} \rightarrow \left\{ \begin{array}{ll} +i & \mbox{ if $\epsilon_{i}=+1$} \\
-i & \mbox{ if $\epsilon_{i}=-1$}
\end{array}\right.
\end{eqnarray*}
For example, if $u=x_1^{-1}x_4x_2x_3^2x_7\in V_7$ then $u$ is
represented by the string \verb=-1 4 2 3 3 7=. In this context the
length of $u$ is equal to the number of integers in its string
representation. We define a
\defterm{chromosome} to be the GA representation of a pair
$(\chi,\zeta)$ of words, and note that each word is naturally of
variable length. Moreover a \defterm{population} is a multiset of a
fixed number $p$ of chromosomes. The GA has two populations in
memory, the
\defterm{current population} and the \defterm{next generation}. As with traditional
GAs the current population contains the chromosomes under
consideration at the current iteration of the GA, and the next
generation has chromosomes deposited into it by the GA which form
the current population on the next iteration. A
\defterm{subpopulation} is a submultiset of a given population.

We use the natural representation for ease of algebraic operation,
acknowledging that faster or more sophisticated data structures
exist, for example the stack-based data structure of
\cite{Wrathall2}. However we believe the simplicity of our
representation yields relatively uncomplicated reproductive
algorithms. In contrast, we believe a stack-based data structure
yields reproductive methods of considerable complexity. We give our
reproductive methods in the next subsection.

Besides normal form reduction of a word $u$ we use
\defterm{pseudo-reduction} of $u$. Let
$\{\,x_{i_{j_1}},x_{i_{j_1}}^{-1},\ldots,x_{i_{j_m}},x_{i_{j_m}}^{-1}\}$
be the generators which would be removed from $u$ if we were to
reduce $u$ to normal form. Pseudo-reduction of $u$ is defined as
simply removing the above generators from $u$. There is no
reordering of the resulting word (as with normal form). For
example, if $u=x_6x_8x_1^{-1}\!x_2x_8^{-1}\!x_2^{-1}\!x_6x_4x_5$
then its \defterm{pseudo-normal form} is
$\tilde{u}=x_6x_1^{-1}\!x_6x_4x_5$ and the normal form of $u$ is
$\overline{u}=x_1^{-1}\!x_4x_6^2x_5$. Clearly, we have
$l(\tilde{u})=l(\overline{u})$. This form is efficiently
computable, with complexity at most that of the algorithm used to
compute the normal form $\overline{u}$. Note, a word is not
assumed to be in any given form unless otherwise stated.

\subsection{Reproduction}\label{reps}
The following reproduction methods are adaptations of standard GA
reproduction methods. The methods act on a subpopulation to give a
child chromosome, which we insert into the next population (more
details are given in section \ref{GA details}).

\begin{enumerate} \item  Sexual
(\defterm{crossover}): by some selection function, input two
parent chromosomes $c_1$ and $c_2$ from the current population.
Choose one random segment from $c_1$, one from $c_2$ and output
the concatenation of the segments. \item Asexual: input a parent
chromosome $c$, given by a selection function, from the current
population. Output one child chromosome by one of the following:
\begin{enumerate}
\item  \defterm{Insertion} of a random generator into a random position
of $c$.
\item  \defterm{Deletion} of a generator at a random position of $c$.
\item  \defterm{Substitution} of a generator located at a random
position in $c$ with a random generator.
\end{enumerate}
\item  Continuance: return several chromosomes $c_1,c_2,\ldots,c_m$ chosen by some
selection algorithm, such that the first one returned is the
``fittest'' chromosome (see the next subsection). This method is
known as \defterm{partially elitist}.
\item  Non-Local Admission: return a random chromosome by
some algorithm.
\end{enumerate}
With the exception of continuance, the methods are repeated for each
child chromosome required.

\subsection{The Cost Function}
In a sense, a cost function induces a partial metric over the search
space to give a measure of the ``distance'' of a chromosome from a
solution. Denote the solution of an instance of the DCSP in section
\ref{prs} by $(x,y)$ and a chromosome by $(\chi,\zeta)$. Let
$E(\chi,\zeta)=\chi a\zeta b^{-1}$; for simplicity we denote this
expression by $E$. The normal form of the above expression is
denoted $\overline{E}$. When $(\chi,\zeta)$ is a solution to an
instance, we have $\overline{E}=\varepsilon$ (the empty word) with
defined length $l(\overline{E})=0$.

The cost function we use is as follows: given a chromosome
$(\chi,\zeta)$ its cost is given by the formula
$C((\chi,\zeta))=l(\overline{E})$. This value is computed for every
chromosome in the current population at each iteration of the GA.
This means we seek to minimise the value of $C((\chi,\zeta))$ as we
iterate the GA.

\subsection{Selection Algorithms}
We realise continuance by roulette wheel selection. This is cost
proportionate. As we will see in Algorithm \ref{GA alg}, we
implicitly require the population to be ordered best cost first.
To this end, write the population as a list
$\{(\chi_1,\zeta_1),\ldots,(\chi_p,\zeta_p)\}$ where
$C(\chi_1,\zeta_1)\leq C(\chi_2,\zeta_2)\leq\ldots\leq
C(\chi_p,\zeta_p)$. Then the algorithm is as follows:

\begin{algorithm}[Roulette Wheel Selection]
\
\begin{itemize}
\item[\textsc{Input:}] The population size $p$; the population
chromosomes $(\chi_i,\zeta_i)$; their costs

$C((\chi_i,\zeta_i))$; and $n_s$, the number of chromosomes to
select
\item[\textsc{Output:}] $n_s$ chromosomes from the population
\end{itemize}
\begin{enumerate}
\item[1.] Let $W\leftarrow\sum_{i=1}^{p}C((\chi_i,\zeta_i))$;
\item[2.] Compute the sequence $\{p_s\}$ such that $p_s((\chi_i,\zeta_i))\leftarrow\frac{C((\chi_i,\zeta_i))}{W}$;
\item[3.] Reverse the sequence $\{p_s\}$;
\item[4.] For $j=1,\ldots,p$, compute $q_j\leftarrow\sum_{i=1}^{j}p_s((\chi_i,\zeta_i))$;
\item[5.] For $t=1,\ldots,n_s$, do

\begin{enumerate}
\item  If $t=1$ output $(\chi_1,\zeta_1)$, the chromosome with least cost. End.
\item Else
\begin{enumerate}
\item  Choose a random $r\in\lbrack0,1\rbrack$;
\item  Output
$(\chi_k,\zeta_k)$ such that $q_{k-1}<r<q_k$. End.
\end{enumerate}
\end{enumerate}
\end{enumerate}
\end{algorithm}

The algorithm respects the requirement that chromosomes with least
cost are selected more often. For crossover we use
\defterm{tournament selection}, where we input three randomly chosen
chromosomes in the current population and select the two with least
cost. If all three have identical cost, then select the first two
chosen. Selection of chromosomes for asexual reproduction is at
random from the current population.

\section{Traceback} In many ways, cost functions
are a large part of a GA. But the reproduction methods often specify
that a random generator is chosen, so reducing the number of
possible choices of generator may serve to guide the GA. We give a
possible approach to reducing this number and term it
\defterm{traceback}. In brief, we take the problem instance given by
the pair $(a,b)$ and use $a$ and $b$ to determine properties of a
feasible solution $(x,y)\in\mathcal{M}((a,b))$ to the instance. This
approach exploits the ``geometry'' of the search space by tracking
the process of reduction of $E$ to its normal form in $V_n$ and
proceeds as follows:

Recall $Y$ and $Z$ respectively denote the set of generators of the
parabolic subgroups $G(Y)$ and $G(Z)$. Suppose we have a chromosome
$(\chi,\zeta)$ at some stage of the GA computation. Form the
expression $E=\chi a\zeta b^{-1}$ associated to the given instance
of the DCSP and label each generator from $\chi$ and $\zeta$ with
its position in the product $\chi\zeta$. Then reduce $E$ to its
normal form $\overline{E}$; during reduction the labels travel with
their associated generators. As a result some generators from $\chi$
or $\zeta$ may be cancelled or not, and the set of labels of the
non-cancelled generators of $\chi$ and $\zeta$ give the original
positions.

The generators in $V_n$ which commute mean that the chromosome may
be split into \defterm{blocks} $\{\beta_i\}$. Each block is formed
from at least one consecutive generator of $\chi$ and $\zeta$ which
move together under reduction of $E$. Let $B$ be the set of all
blocks from the above process. Now a block $\beta_m\in B$ and a
position $q$ (which we call the \defterm{recommended position}) at
either the left or right end of that block are randomly chosen.
Depending upon the position chosen, take the subword $\delta$
between either the current and next block $\beta_{m+1}$ or the
current and prior block $\beta_{m-1}$ (if available). If there is
just one block, then take $\delta$ to be between $\beta_1$ and the
end or beginning of $\overline{E}$.

Then identify the word $\chi$ or $\zeta$ from which the position $q$
originated and its associated generating set $S=Y$ or $S=Z$. The
position $q$ is at either the left or right end of the chosen block.
So depending on the end of the block chosen, randomly select the
inverse of a generator from $RF(\delta)\cap S$ or $FL(\delta)\cap
S$. Call this the \defterm{recommended generator} $g$. Note if both
$\chi$ and $\zeta$ are entirely cancelled (and so $B$ is empty), we
return a random recommended generator and position.

With these, the insertion algorithm inserts the inverse of the
generator on the appropriate side of the recommended position in
$\chi$ or $\zeta$. In the cases of substitution and deletion, we
substitute the recommended generator or delete the generator at the
recommended position. We now give an example for the DCSP on
$V_{10}$ with the two parabolic subgroups of $V(Y)=V_7$ and
$V(Z)=V_{10}$.

\subsubsection*{Example of Traceback on a Given Instance}
Take the short DCSP instance
\[
(a,b)=(x_2^2x_3x_4x_5x_4^{-1}x_7x_6^{-1}x_9x_{10},\,\,x_2^2x_4x_5x_4^{-1}x_3x_7x_6^{-1}x_{10}x_9)
\]
and let the current chromosome be
$(\chi,\zeta)=(x_3x_2^{-1}x_3^{-1}x_5x_7,\,\,x_5x_2x_3x_7^{-1}x_{10})$.
Represent the labels of the positions of the generators in $\chi$
and $\zeta$ by the following numbers immediately above each
generator:
\begin{eqnarray*}
\begin{array}{ccccc|ccccc}
  0 & 1 & 2 & 3 & 4 \,\,&\,\, 5 & 6 & 7 & 8 & 9 \\
  x_3 & x_2^{-1} & x_3^{-1}& x_5 & x_7 \,\, & \,\, x_5 & x_2 & x_3& x_7^{-1}&
  x_{10}\\
\end{array}
\end{eqnarray*}
Forming $E$ and reducing it to its Knuth-Bendix normal form gives
\begin{eqnarray*}
\overline{E} &=&
\begin{array}{ccccccccccccccc}
  0 & 1 & 2 & & & & & 3 & & & & 4 & \\
  x_3 & x_2^{-1} & x_3^{-1} & x_2 & x_2 & x_3 & x_2^{-1} & x_5 & x_4 & x_5 & x_4^{-1} & x_7 & x_7\\
  & 5 & & 8 & & & & & & & 9 & & \\
  x_6^{-1} & x_5 & x_4 & x_7^{-1} & x_6 & x_5^{-1} & x_4^{-1} & x_7^{-1} & x_9 & x_{10} & x_{10} & x_9^{-1} & x_{10}^{-1} \\
\end{array}
\end{eqnarray*}
which contains eight remaining generators from $(\chi,\zeta)$. Take
cost to be $ C((\chi,\zeta))=l(\overline{E})=26$, the number of
generators in $\overline{E}$ above. There are three blocks for
$\chi$:
\[
\beta _{1}=
\begin{array}{ccc}
0 & 1 & 2 \\
x_{3} & x_{2}^{-1} & x_{3}
\end{array}
,\,\,\beta _{2}=
\begin{array}{c}
3 \\
x_{5}
\end{array}
,\,\,\beta _{3}=
\begin{array}{c}
4 \\
x_{7}
\end{array}
\]
and three for $\zeta$:
\[
\beta _{4}=
\begin{array}{c}
5 \\
x_{5}
\end{array}
,\,\,\beta _{5}=
\begin{array}{c}
8 \\
x_{7}^{-1}
\end{array}
,\,\,\beta _{6}=
\begin{array}{c}
9 \\
x_{10}
\end{array}
\]
Suppose we choose position eight, which is in $\zeta$ and is block
$\beta_5$. This is a block of length one; we may take the word to
the left or the right as our choice for $\delta$.

Suppose we choose the word to the right, so
$\delta=x_6x_5^{-1}x_4^{-1}x_7^{-1}x_9x_{10}$ and in this case,
$S=\{x_1,\ldots,x_{10}\}$. So we choose a random generator from
$FL(\delta)\cap S=\{x_6,x_9\}$. Choose $g=x_6^{-1}$ and so $\zeta$
becomes $\zeta'=x_5x_2x_3x_7^{-1}x_6^{-1}x_{10}$, with $\chi'=\chi$.
The cost becomes $C((\chi',\zeta'))=l(\overline{\chi' a\zeta'
b^{-1}})=25$. Note that we could have taken any block and the
permitted directions to create $\delta$. In this case, there are
eleven choices of $\delta$, clearly considerably fewer than the
total number of subwords of $\overline{E}$. Traceback provides a
significant increase in performance over merely random selection
(this is easily calculated in the above example to be by a factor of
$38$).

\section{Setup of the Genetic Algorithm}\label{GA details}

\subsection{Specification of Output Alphabet}
Let $n=2m$ for some integer $m>1$. Define the subsets of generators
$Y=\{x_1,\ldots,x_{m-1}\},\,\,Z=\{x_{m+2},\ldots,x_n\}$ and two
corresponding parabolic subgroups
$G(Y)=\left<Y\right>,G(Z)=\left<Z\right>$. Clearly $G(Y)$ and $G(Z)$
commute as groups: if we take any $m>1$ and any words $x_y\in G(Y)$,
$x_z\in G(Z)$ then $x_yx_z=x_zx_y$. We direct the interested reader
to \cite{KoSECANTS02} for information on the importance of the
preceding statement. Given an instance $(a,b)$ of the DCSP with
parabolic subgroups as above, we will seek a representative for each
of the two words $x\in G(Y)$ and $y\in G(Z)$ that are a solution to
the DCSP. Let us label this problem $(P)$.

\subsection{The Algorithm and its Parameters}
Given a chromosome $(\chi,\zeta)$ we choose crossover to act on
either $\chi$ or $\zeta$ at random, and fix the other component of
the chromosome. Insertion is performed according to the position in
$\chi$ or $\zeta$ given by traceback and substitution is with a
random generator, both such that if the generator chosen cancels
with a neighbouring generator from the word then another random
generator is chosen. We choose to use pseudo-normal form for all
chromosomes to remove all redundant generators while preserving the
internal ordering of $(\chi,\zeta)$.

By experiment, GA behaviour and performance is mostly controlled by
the
\defterm{parameter set} chosen. A parameter set is specified by the
population size $p$ and numbers of children begat by each
reproduction algorithm. The collection of numbers of children is
given by a multiset of non-negative integers $P=\{p_i\}$, where
$\sum p_i=p$ and each $p_i$ is given, in order, by the number of
crossovers, substitutions, deletions, insertions, selections and
random chromosomes. The GA is summarised by the following algorithm:

\begin{algorithm}[GA for DCSP]\label{GA alg}
\
\begin{itemize}
\item[\textsc{Input:}] The parameter set, words $a,b$ and their
lengths $l(a),l(b)$, suicide control $\sigma$, initial length $L_I$
\item[\textsc{Output:}] A solution $(\chi,\zeta)$ or timeout; $i$, the
number of populations
\end{itemize}

\begin{enumerate}
\item Generate the initial population $P_0$, consisting of $p$
random (unreduced) chromosomes $(\chi,\zeta)$ of initial length
$L_I$;

\item $i\leftarrow 0$;

\item Reduce every chromosome in the population to its
pseudo-normal form.

\item While $i<\sigma$ do

\begin{enumerate}
\item For $j=1,\ldots,p$ do

\begin{enumerate}
\item Reduce each pair $(\chi_j,\zeta_j)\in P_i$ to its
pseudo-normal form $(\tilde{\chi_j},\tilde{\zeta_j})$;

\item Form the expression $E=\tilde{\chi_j}\,a\,\tilde{\zeta_j}\,b^{-1}$;

\item Perform the traceback algorithm to give $C((\chi_j,\zeta_j))$,
recommended generator $g$ and recommended position $q$;
\end{enumerate}
\item Sort current population $P_i$ into least-cost-first order
and label the chromosomes
$(\tilde{\chi_1},\tilde{\zeta_1}),\ldots,(\tilde{\chi_p},\tilde{\zeta_p})$;

\item If the cost of $(\tilde{\chi_1},\tilde{\zeta_1})$ is zero then
return solution $(\chi_1,\zeta_1)$ and \textsc{END.}

\item $P_{i+1}\leftarrow\emptyset$;

\item For $j=1,\ldots,p$ do
\begin{enumerate}
\item Using the data obtained in step $4(a)(iii)$, perform the appropriate reproductive algorithm on
$(\tilde{\chi_j},\tilde{\zeta_j})$ and denote the result
$(\chi_j',\zeta_j')$;

\item $P_{i+1}\leftarrow P_{i+1}\cup\{(\chi_j',\zeta_j')\}$;
\end{enumerate}
\item $i\leftarrow i+1$.

\end{enumerate}
\item Return failure. \textsc{END.}
\end{enumerate}
\end{algorithm}
The positive integer $\sigma$ is an example of a
\defterm{suicide control}, where the GA stops (suicide) if more than
$\sigma$ populations have been generated. In all cases here,
$\sigma$ is chosen by experimentation; GA runs that continued beyond
$\sigma$ populations were unlikely to produce a successful
conclusion. By deterministic search we found a population size of
$p=200$ and parameter set $P=\{5,33,4,128,30,0\}$ for which the GA
performs well when $n=10$. We observed that the GA exhibits the
well-known common characteristic of sensitivity to changes in
parameter set; we consider this in future work. We found an optimal
length of one for each word in our initial population, and now
devote the remainder of the paper to our results of testing the GA
and analysis of the data collected.

\subsection{Method of Testing} We wished to test the performance of
the GA on ``randomly generated'' instances of problem $(P)$. Define
the length of an instance of $(P)$ to be the set of lengths
$\{l(\overline{a}),l(\overline{x}),l(\overline{y})\}$ of words
$a,x,y\in V_n$ used to create that instance. Each of the words $a,x$
and $y$ are generated by simple random walk on $V_n$. To generate a
word $\overline{u}$ of given length $k=l(\overline{u})$ firstly
generate the unreduced word $u_1$ with unreduced length $l(u_1)=k$.
Then if $l(\overline{u_1})<k$, generate $u_2$ of unreduced length
$k-\l(\overline{u_1})$, take $u_1u_2$ and repeat this procedure
until we produce a word $u=u_1u_2\ldots\ u_r$ with $l(\overline{u})$
equal to the required length $k$.

We identified two key input data for the GA: the length of an
instance of $(P)$ and the group rank, $n$. Two types of tests were
performed, varying these data:
\begin{enumerate}
\item Test of the GA with long instances while keeping the rank small;

\item Test of the GA with instances of moderate length while increasing the
rank.
\end{enumerate}
The algorithms and tests were developed and conducted in GNU C++ on
a Pentium IV 2.53\,GHz computer with 1GB of RAM running Debian Linux
3.0.

\subsection{Results} Define the
\defterm{generation count} to be the number of populations (and so
iterations) required to solve a given instance; see the counter
$i$ in Algorithm \ref{GA alg}. We present the results of the tests
and follow this in section \ref{DC} with discussion of the
results.

\subsubsection*{Increasing Length} We tested the GA on eight
randomly generated instances (I1)--(I8) with the rank of $V_n$ set
at $n=10$. The instances (I1)--(I8) were generated beginning with
$l(\overline{a})=128$ and $l(\overline{x})=l(\overline{y})=16$ for
instance (I1) and progressing to the following instance by doubling
the length $l(\overline{a})$ or both of the lengths
$l(\overline{x})$ and $l(\overline{y})$. The GA was run ten times on
each instance and the mean runtime $\overline{t}$ in seconds and
mean generation count $\overline{g}$ across all runs of that
instance was taken. For each collection of runs of an instance we
took the standard deviation $\sigma_g$ of the generation counts and
the mean time in seconds taken to compute each population. A summary
of results is given in Table \ref{fig2}.
\begin{table}
\centering \caption{Results of increasing instance lengths for
constant rank $n=10$.}
\begin{tabular}{c|ccc|cc|cc}
\hline
Instance & $l(\overline{a})$ & $l(\overline{x})$ & $l(\overline{y})$ & $\overline{g}$ & $\overline{t}$ & $\sigma_g$ & sec/gen\\
\hline\hline
I1 & 128 & 16 & 16 & 183 & 59 & 68.3 & 0.323\\
I2 & 128 & 32 & 32 & 313 & 105 & 198.5 & 0.339\\
I3 & 256 & 64 & 64 & 780 & 380 & 325.5 & 0.515\\
I4 & 512 & 64 & 64 & 623 & 376 & 205.8 & 0.607\\
I5 & 512 & 128 & 128 & 731 & 562 & 84.4 & 0.769\\
I6 & 1024 & 128 & 128 & 1342 & 801 & 307.1 & 0.598\\
I7 & 1024 & 256 & 256 & 5947 & 5921 & 1525.3 & 1.004\\
I8 & 2048 & 512 & 512 & 14805 & 58444 & 3576.4 & 3.849\\
\hline
\end{tabular}
\label{fig2}
\end{table}

\noindent\textbf{Increasing Rank} These tests were designed to keep
the lengths of computed words relatively small while allowing the
rank $n$ to increase. We no longer impose the condition of
$l(\overline{x})=l(\overline{y})$. Take $s$ to be the arithmetic
mean of the lengths of $\overline{x}$ and $\overline{y}$. Instances
were constructed by taking $n=10,20$ or $40$ and generating random
$a$ of maximal length 750, random $x$ and $y$ of maximal length 150
and then reducing the new $b=xay$ to its normal form $\overline{b}$.

We then ran the GA once on each of $505$ randomly generated
instances for $n=10$, with $145$ instances for $n=20$ and $52$
instances for $n=40$. We took the time $t$ in seconds to produce a
solution and the respective generation count $g$. The data collected
is summarised on Table \ref{fig1} by grouping the length $s$ of
instance into intervals of length fifteen. For example, the range
75--90 means all instances where $s\in[75,90)$. Across each interval
we computed the means $\overline{g}$ and $\overline{t}$ along with
the standard deviation $\sigma_g$. We now give a brief discussion of
the results and some conjectures, and then conclude our work.

\begin{table}
\centering
\caption{Results of increasing rank from $n=10$ (upper
rows) to $n=20$ (centre rows) and $n=40$ (lower rows).}
\begin{tabular}{|c|c|c|c|c|c|c|c|c|c|}
\hline $s$ & 15--30 & 30--45 & 45--60 & 60--75 & 75--90 & 90--105 & 105--120 & 120--135 & 135--150 \\
\hline\hline
$\overline{g}$ & 227 & 467 & 619 & 965 & 1120 & 1740 & 1673 & 2057 & 2412 \\
$\overline{t}$ & 44 & 94 & 123 & 207 & 244 & 384 & 399 & 525 & 652 \\
\hline
$\overline{g}$ & 646 & 2391 & 2593 & 4349 & 4351 & 8585 & 8178 & 8103 & 10351 \\
$\overline{t}$ & 251 & 897 & 876 & 1943 & 1737 & 3339 & 3265 & 4104 & 4337 \\
\hline
$\overline{g}$ & 1341 & 1496 & 2252 & 1721 & 6832 & 14333 & 14363 & - & - \\
$\overline{t}$ & 949 & 1053 & 836 & 1142 & 5727 & 10037 & 11031 & - & - \\
\hline
\end{tabular}
\label{fig1}
\end{table}

\subsection{Discussion and Conclusion}\label{DC} Firstly, the mean
times given on Tables \ref{fig2} and \ref{fig1} depend upon the time
complexity of the underlying algebraic operations. We conjecture for
$n=10$ that these have time complexity no greater than $O(k \log k)$
where $k$ is the mean length of all words across the entire run of
the GA that we wish to reduce.

Table \ref{fig2} shows we have a good method for solving large scale
problems when the rank is $n=10$. By Table \ref{fig1} we observe the
GA operates very well in most cases across problems where the mean
length of $x$ and $y$ is less than $150$ and rank at most forty.
Fixing $s$ in a given range, the mean generation count increases at
an approximately linearithmic rate as $n$ increases. This seems to
hold for all $n$ up to forty, so we conjecture that for a mean
instance of problem $(P)$ with given rank $n$ and instance length
$s$ the generation count for an average run of the GA lies between
$O(sn)$ and $O(sn\log n)$. This conjecture means the GA generation
count depends linearly on $s$ (for brevity, we omit the statistical
evidence here).

As $n$ increases across the full range of instances of $(P)$,
increasing numbers of suicides tend to occur as the GA encounters
increasing numbers of local minima. These may be partially explained
by observing traceback. For $n$ large, we are likely to have many
more blocks than for $n$ small (as the likelihood of two arbitrary
generators commuting is larger). While traceback is much more
efficient than a purely random method, there are more chances to
read $\delta$ between blocks. Indeed, there may be so many possible
$\delta$ that it takes many GA iterations to reduce cost. By
experience of this situation, non-asexual methods of reproduction
bring the GA out of some local minima. Consider the following
typical GA output, where the best chromosomes from populations $44$
and $64$ (before and after a local minimum) are:
\begin{verbatim}
Gen 44 (c = 302) : x = 9 6 5 6 7 4 5 -6 7 5 -3 -3 (l = 12)

y = -20 14 12 14 -20 -20 (l = 6)
\end{verbatim}
\begin{verbatim}Gen 64 (c = 300) : x = 9 8 1 7 6 5 6 7 4 5 -6 7 9 5 -3 -3 (l = 16)

y = 14 12 12 -20 14 15 -14 -14 -16 17 15 14 -20 15 -19 -20 -20 -19
-20 18 -17 -16 (l = 22)
\end{verbatim}
In this case, cost reduction is not made by a small change in
chromosome length, but by a large one. We observe that the cost
reduction is made when a chromosome from lower in the ordered
population is selected and then mutated, as the new chromosome at
population $64$ is far longer. In this case it seems traceback acts
as a topological sorting method on the generators of the equation
$E$, giving complex systems of cancellation in $E$ which result in a
cost deduction greater than one. This suggests that finetuning the
parameter set to focus more on reproduction lower in the population
and reproduction which causes larger changes in word length may
improve performance. Indeed, \cite{Bremer62} conjectures that
\begin{quote}
``It seems plausible to conjecture that sexual mating has the
purpose to overcome situations where asexual evolution is
stagnant.'' \flushright{\textbf{Bremermann \cite[p.~102]{Bremer62}}}
\end{quote}
This implies the GA performs well in comparison to asexual
hillclimbing methods. Indeed, this is the case in practice: by
making appropriate parameter choices we may simulate such a
hillclimb, which experimentally encounters many more local minima.
These local minima seem to require substantial changes in the form
of $\chi$ and $\zeta$ (as above); this cannot be done by mere
asexual reproduction.

Meanwhile, coupled with a concept of ``growing'' solutions, we have
at least for reasonable values of $n$ an indication of a good
underlying deterministic algorithm based on traceback. Indeed, such
deterministic algorithms were developed in \cite{BEKR03} as the
result of analysis of experimental data in our work. This hints that
the search space has a ``good'' structure and may be exploited by
appropriately sensitive GAs and other artificial intelligence
technologies in our framework.


\begin{thebibliography}{XX}
\bibitem{BBB03} R. F. Booth, D. Y. Bormotov, A. V. Borovik,
Genetic Algorithms and Equations in Free Groups and Semigroups,
Contemp. Math. \textbf{349} (2004), 63--80.

\bibitem{BEKR03} A. V. Borovik, E. S. Esyp, I. V. Kazatchkov, V. N.
Remeslennikov, Divisibility Theory and Complexity of Algorithms for
Free Partially Commutative Groups, Contemp. Math. \textbf{378}
(Groups, Languages, Algorithms), 2005.

\bibitem{Bremer62} H. J. Bremermann, Optimization Through
Evolution and Recombination, Self-Organizing Systems (M. C. Yovits
et al., eds.), Washington, Spartan Books (1962), 93--106.

\bibitem{Holland1975} J. Holland, Adaptation in Natural and
Artificial Systems (5th printing), MIT Press, Cambridge,
Massachusetts, 1998.

\bibitem{KoSECANTS02} K. -H. Ko, Braid Group and Cryptography,
19th SECANTS, Oxford, 2002.

\bibitem{KnuthBendix} D. Knuth, P. Bendix, Simple Word Problems in
Universal Algebra, Computational Problems in Abstract Algebras (J.
Leech, ed.), Pergamon Press 1970, 263--297.

\bibitem{Miasnikov99} A. D. Miasnikov, Genetic Algorithms and the
Andrews-Curtis Conjecture, Internat. J. Algebra Comput. \textbf{9}
(1999), no. 6, 671--686.

\bibitem{2Miasnikov} A. D. Miasnikov, A. G. Myasnikov, Whitehead
Method and Genetic Algorithms, Contemp. Math. \textbf{349} (2004),
89--114.

\bibitem{Michalewicz} Z. Michalewicz, Genetic Algorithms + Data
Structures = Evolution Programs (3rd rev. and extended ed.),
Springer-Verlag, Berlin, 1996.

\bibitem{VanWyk} L. VanWyk, Graph Groups are Biautomatic, J. Pure Appl. Algebra
\textbf{94} (1994), no. 3, 341--352.

\bibitem{Vershik} A. Vershik, S. Nechaev, R. Bikbov, Statistical
Properties of Braid Groups in Locally Free Approximation, Comm.
Math. Phys. \textbf{212} (2000), 59--128.

\bibitem{Wrathall1} C. Wrathall, The Word Problem for Free
Partially Commutative Groups, J. Symbolic Comp. \textbf{6} (1988),
99--104.

\bibitem{Wrathall2} C. Wrathall, Free partially commutative
groups, Combinatorics, Computing and Complexity (Tianjing and
Beijing, 1988) 195--216, Math. Appl (Chin. Ser. 1) Kluwer Acad.
Publ., Dordrecht, 1989.
\end{thebibliography}
\end{document}